\documentclass[12pt, draft]{article}
\usepackage{amsmath, amsfonts, amsthm, amssymb}
\newtheorem{theorem}{Theorem}[section]
\newtheorem{prove}{Proof of Theorem}[section]
\newtheorem{lemma}[theorem]{Lemma}
\newtheorem{definition}{Definition}[section]
\newtheorem{proposition}[definition]{Proposition}
\begin{document}
\title{\bf ON A SUBCLASS OF 5-DIMENSIONAL SOLVABLE  LIE
ALGEBRAS WHICH HAVE 3-DIMENSIONAL COMMUTATIVE DERIVED IDEAL}
\author{\bf Le Anh Vu\\Department of Mathematics and Informatics\\
Ho Chi Minh City University of Pedagogy\\
e-mail:leanhvu@hcmup.edu.vn} \footnotetext{

2000AMS Mathematics Subject Classification: Primary 22E45, Secondary
46E25, 20C20.}
\date{}
\maketitle
\vskip 1cm
\begin{abstract}
    The paper presents a subclass  of the class of MD5-algebras and
MD5-groups,  i.e.,  five dimensional solvable Lie algebras and Lie
groups such that their orbits in the co-adjoint representation
(K-orbit) are orbit of zero or maximal dimension. The main results
of the paper is the classification up to an isomorphism  of all
MD5-algebras $\mathcal{G}$ with the derived ideal ${\mathcal{G}}^{1}
:= [\mathcal{G}, \mathcal{G}]$ is a 3-dimensional commutative Lie
algebra. \end{abstract}
\subsection*{Introduction}
    In 1962, studying theory of representations,
Kirillov [3] introduced the Orbit Method. This method quickly became
the most important method in the theory of representations of Lie
groups. The Kirillov's Orbit Method immediately was expanded by
Kostant, Auslander, Do Ngoc Diep, etc. Using the Kirillov's Orbit
Method, we can obtain all the unitary irreducible representations of
solvable and simply connected Lie groups. The importance of
Kirillov's Orbits Method is the co-adjoint representation
(K-representation). Therefore, it is meaningful to study the
K-representation in the theory of  representations of Lie groups.

The structure of solvable Lie groups and Lie algebras is not to
complicated, although the complete classification of them is
unresolved up to now. In 1980,  studying the Kirillov's Orbit
Method, D. N. Diep [2] introduced the class of  Lie groups and Lie
algebras MD. Let G be an n-dimensional real Lie group. It is called
an MDn-group iff its orbits in the K-representation (i.e. K-orbits)
are orbits of dimension zero or maximal dimension. The corresponding
Lie algebra of G is called an MDn-algebra. Thus, classification and
study of K-representation of the class of MDn-groups and
MDn-algebras is the problem of great interest. Because of all Lie
algebras of n dimension (with $n \leq 3$) were listed easily, we
have to consider MDn-groups and MDn-algebras with $n \geq 4$.

In 1984, Dao Van Tra ([5]) was listed all MD4-algebras. In 1992, all
MD4-algebras were classified up to an isomorphism by the author (see
[6], [7], [8]). Until now, no complete classification of
MDn-algebras with $n \geq 5$ is known. Three first examples of
MD5-algebras and MD5-groups can be found in [9] and some different
MD5-algebras and MD5-groups can be found in [10]. In this paper we
shall give the classification up to an isomorphism of all
MD5-algebras $\mathcal{G}$ with the derived ideal ${\mathcal{G}}^{1}
:= [\mathcal{G}, \mathcal{G}]$ is a 3-dimensional commutative Lie
algebra. The complete classification of all MD5-algebras will be
presented in the next paper.

\section{Preliminaries}
At first, we recall in this section some preliminary results and
notations which will be used later. For details we refer the reader
to References [2], [3], [4].

\subsection{Lie Groups and Lie Algebras}
\begin{definition}  A real Lie group of dimension n is a $C^{\infty}$-manifold G endowed
with a group structure such that the map $(g, h) \mapsto g.h^{-1}$
from G$\times$G into G is $C^{\infty}$-differentiable.
\end{definition}
\begin{definition}  A real Lie algebra $\mathcal{G}$ of dimension n is an
n-dimensional real vector space together with a skew-symmetric
bilinear map (X, Y) $\mapsto$ [X, Y] from ${\mathcal{G}}\times
{\mathcal{G}}$ into $\mathcal{G}$ (which is called the Lie bracket)
such that the following Jacobi identity is satisfied : [[X,Y],Z] +
[[Y,Z],X] + [[Z,X],Y] = 0 for every X, Y, Z $\in \mathcal{G}$.
\end{definition}
\subsection{The co-adjoint Representation, K-orbits \\ MDn-Groups and MDn-Algebras}
    Each Lie group G determines a Lie algebra Lie(G) = $\mathcal{G}$ as the tangent
space $ T_{e}G $ of G at the identity with the Lie bracket is
defined by the commutator. Inversely, each real Lie algebra
$\mathcal{G}$ is associated to one connected and simply connected
Lie group G such Lie(G) = $\mathcal{G}$. For each g $\in$ G, we
denote the internal automorphism associated to g by $A_{(g)}$. So $
A_{(g)}: G \longrightarrow G $ is defined as follows
$$A_{(g)}(x):=\, g.x.g^{-1},\, \forall x \in G.$$
This automorphism induces the following map
$${A_{(g)}}_{*}:\mathcal{G}\longrightarrow \mathcal{G}\qquad \qquad \qquad \qquad $$
$$\qquad\qquad \qquad \qquad\qquad X\longmapsto {A_{(g)}}_{*}(X):\, =\;\frac{d}{dt}
[g.exp(tX)g^{-1}]\mid_{t=0}$$ which is called the tangent map of
$A_{(g)}$.

\begin{definition}The action
$$Ad :G\longrightarrow Aut(\mathcal{G})$$
$$\qquad\qquad\qquad g\longmapsto Ad(g):\, =\;  {A_{(g)}}_{*}$$
is called the adjoint representation of G in $\mathcal{G}$.
\end{definition}

\begin{definition}The action
$$K:G\longrightarrow Aut(\mathcal{G}^{*})$$
$$g\longmapsto K_{(g)}$$
such that
$$\langle K_{(g)}F,X\rangle :\,=\langle F, Ad(g^{-1})X\rangle ;
\quad(F\in {\mathcal{G}}^{*},\, X \in \mathcal{G})$$ is called the
co-adjoint representation of G in $\mathcal{{G}^{*}}$.
\end{definition}
\begin{definition} Each orbit of the co-adjoint representation
of G is called a K-orbit of G. \end{definition} Thus, for every $F
\in \mathcal{G}^{*}$, the K-orbit containing $F$ is defined as
follows
$${\varOmega}_{F}:= \{K_{(g)}F / g \in G \}.$$  The dimension of
every K-orbit of G is always even. In order to define the dimension
of the K-orbits ${\varOmega}_{F}$, it is useful to consider the
skew-symmetric bilinear form $B_{F}$ on $\mathcal{G}$ as follows
$$ B_{F}(X, Y) := \langle F, [X, Y]\rangle; \, \forall\, X, Y \in
\mathcal{G}.$$ Denote the stabilizer of $F$ under the co-adjoint
representation of G in $\mathcal{{G}^{*}}$ by $G_{F}$ and
${\mathcal{G}}_{F}$ := Lie($G_{F}$). We shall need in the sequel the
following assertion.
\begin{proposition}[see {[3]}] Ker$B_{F}$ = ${\mathcal{G}}_{F}$ and
dim${\varOmega}_{F}$ = dim$\mathcal{G} - dim{\mathcal{G}}_{F}$.
\hfill{$\square$} \end{proposition}
\begin{definition}[{see [2]}]An MDn-group is an n-dimensional real solvable
Lie group such that its K-orbits are orbits of dimension zero or
maximal dimension. The Lie algebra of an MDn-group is called an
MDn-algebra. \end{definition} The following proposition give a
necessary condition in order that a Lie algebra belongs to the class
of all MD-algebras.
\begin{proposition}[{see [4]}] Let $\mathcal{G}$ be an MD-algebra. Then its
second derived ideal ${\mathcal{G}}^{2} := [[\mathcal{G},
\mathcal{G}], [\mathcal{G}, \mathcal{G}]]$ is commutative.
\hfill{$\square$} \end{proposition} Note, however, that the converse
of this statement in general is not hold. In other words, the above
necessary condition is not sufficient one.
\section{The Main Result}
From now on, $\mathcal{G}$ will denote an Lie algebra of  dimension
5. We always choose a suitable basis $(\nobreak X_{1}, X_{2}, X_{3},
X_{4}, X_{5} \nobreak)$ in $\mathcal{G}$. Then $\mathcal{G}$
isomorphic to ${\bf R}^{5}$ as a real vector space. The notation
${\mathcal{G}}^{*}$ will mean the dual space of $\mathcal{G}$.
Clearly ${\mathcal{G}}^{*}$ can be identified with ${\bf R}^{5}$ by
fixing in it the basis $(\nobreak X_{1}^{*}, X_{2}^{*}, X_{3}^{*},
X_{4}^{*}, X_{5}^{*}\nobreak )$ dual to the basis $(\nobreak X_{1},
X_{2}, X_{3}, X_{4}, X_{5} \nobreak)$.

\begin{theorem}
Let $\mathcal{G}$ be an MD5-algebra with ${\mathcal{G}}^{1} :=
[\mathcal{G}, \mathcal{G}] \cong {\bf {R}^{3}}$ ( the 3-dimensional
commutative Lie algebra ).
\begin{description}
    \item[I.] Assume that $\mathcal{G}$ is decomposable. Then
    $\mathcal{G} \cong \mathcal{H} \oplus {\bf R}$, where
    $\mathcal{H}$ is an MD4-algebra.
    \item[II.] Assume that $\mathcal{G}$ is indecomposable. Then
    we can choose a suitable basis $( X_{1}, X_{2}, X_{3},
X_{4}, X_{5} )$ of \, $\mathcal{G}$ such that
    ${\mathcal{G}}^{1} = {\bf R}.X_{3} \oplus
{\bf R}.X_{4} \oplus {\bf R}.X_{5} \equiv {\bf {R}^{3}}$,\,
$ad_{X_{1}} = 0 $,\, $ad_{X_{2}} \in End({\mathcal{G}}^{1}) \equiv
    Mat_{3}(\bf R)$;\, $[X_{1}, X_{2}] = X_{3}$ and
    $\mathcal{G}$ is isomorphic to one and only one of the following Lie algebras:
       \begin{description}
         \item[1.]${\mathcal{G}}_{5,3,1({\lambda}_{1}, {\lambda}_{2})}$ :
                    $$ad_{{X}_2} = \begin{pmatrix} {{\lambda}_1}&0&0\\
                    0&{{\lambda}_2}&0\\0&0&1 \end{pmatrix}; \quad
                     {\lambda}_1, {\lambda}_2 \in {\bf R}\setminus \lbrace
                     1\rbrace, \, {\lambda}_1 \neq {\lambda}_2 \neq 0 .$$ \vskip 0.5cm
         \item[2.]${\mathcal{G}}_{5,3,2(\lambda)}$ :
                    $$ad_{{X}_2} = \begin{pmatrix} 1&0&0\\
                    0&1&0\\0&0&{\lambda} \end{pmatrix}; \quad
                    {\lambda} \in {\bf R}\setminus \lbrace 0, 1 \rbrace .$$ \vskip 0.5cm
         \item[3.]${\mathcal{G}}_{5,3,3(\lambda)}$ :
                    $$ad_{{X}_2} = \begin{pmatrix} {\lambda}&0&0\\
                    0&1&0\\0&0&1 \end{pmatrix}; \quad
                    {\lambda} \in {\bf R}\setminus \lbrace 1 \rbrace .$$ \vskip 0.5cm
         \item[4.]${\mathcal{G}}_{5,3,4}$ :
                    $$ad_{{X}_2} = \begin{pmatrix} 1&0&0\\
                    0&1&0\\0&0&1 \end{pmatrix}.$$ \vskip 0.5cm
         \item[5.]${\mathcal{G}}_{5,3,5(\lambda)}$ :
                    $$ad_{{X}_2} = \begin{pmatrix} {\lambda}&0&0\\
                    0&1&1\\0&0&1 \end{pmatrix}; \quad
                    {\lambda} \in {\bf R}\setminus \lbrace 1 \rbrace .$$ \vskip 0.5cm
         \item[6.]${\mathcal{G}}_{5,3,6(\lambda)}$ :
                    $$ad_{{X}_2} = \begin{pmatrix} 1&1&0\\
                    0&1&0\\0&0&{\lambda} \end{pmatrix}; \quad
                    {\lambda} \in {\bf R}\setminus \lbrace 0, 1 \rbrace .$$ \vskip 0.5cm
         \item[7.]${\mathcal{G}}_{5,3,7}$ :
                    $$ad_{{X}_2} = \begin{pmatrix} 1&1&0\\
                    0&1&1\\0&0&1 \end{pmatrix}.$$ \vskip 0.5cm
         \item[8.]${\mathcal{G}}_{5,3,8(\lambda, \varphi)}$ :
                    $$ad_{{X}_2} = \begin{pmatrix} cos{\varphi}&-sin{\varphi}&0\\
                    sin{\varphi}&cos{\varphi}&0\\0&0&\lambda \end{pmatrix}; \quad
                     \lambda \in {\bf R}\setminus \lbrace
                     0\rbrace, \, \varphi \in (0, \pi) .$$ \vskip 0.5cm
       \end{description}
\end{description}
\end{theorem}\vskip3mm
In order to prove Theorem 2.1 we need some lemmas.
\begin{lemma} Under the above notation. We have $ad_{X_{1}}\circ ad_{X_{2}} = ad_{X_{2}}\circ
ad_{X_{1}}$. \end{lemma}
\begin{proof}  Using the Jacobi
identity for $X_{1}, X_{2}$ and $X_{i} ( i = 3, 4, 5 $ ), we have
\begin{align} \notag &\quad [[X_{1}, X_{2}], X_{i}]
+ [[X_{2}, X_{i}], X_{1}] + [[X_{i}, X_{1}], X_{2}] = 0\\
\notag \Leftrightarrow & \quad [X_{1}, [X_{2}, X_{i}]] - [X_{2},
[X_{1}, X_{i}]] = 0\\ \notag \Leftrightarrow & \quad ad_{X_{1}}\circ
ad_{X_{2}} (X_{i}) =
ad_{X_{2}}\circ ad_{X_{1}} (X_{i});\, i = 3, 4, 5 \\
\notag \Leftrightarrow & \quad ad_{X_{1}}\circ ad_{X_{2}} =
ad_{X_{2}}\circ ad_{X_{1}}.  \end{align} \end{proof}

\begin{lemma}[see{[2], [4]}] If \, $\mathcal{G}$\, is an MD-algebra and F
$\in {\mathcal{G}}^{*}$ is not perfectly vanishing on
${\mathcal{G}}^{1}$, i.e. there exists $U \in {\mathcal{G}}^{1}$
such that $\langle F, U \rangle \neq 0$, then the K-orbit
${\Omega}_{F}$ is the one of maximal dimension.
\end{lemma}
\begin{proof} Assume that
${\Omega}_{F}$ is not a K-orbit of maximal dimension, i.e.
$dim{\Omega}_{F} = 0$. This means that
$$dim{\mathcal{G}}_{F} = dim{\mathcal{G}} - dim{\Omega}_{F} =
dim{\mathcal{G}}.$$ So $KerB_{F} = {\mathcal{G}}_{F} = \mathcal{G}
\supset {\mathcal{G}}^{1}$ and F is perfectly vanishing on
${\mathcal{G}}^{1}$. This contradicts the supposition of the lemma.
Therefore ${\Omega}_{F}$ is a K-orbit of maximal
dimension.\end{proof}

We are now in a position to prove the main theorem of the paper.

\begin{prove} \end{prove}  Firstly, we can always choose
some basis $( X_{1}, X_{2}, X_{3}, X_{4}, X_{5} )$ of $\mathcal{G}$
such that ${\mathcal{G}}^{1} = {\bf R}.X_{3} \oplus {\bf R}.X_{4}
\oplus {\bf R}.X_{5} \equiv {\bf {R}^{3}}$;\, $ad_{X_{1}},\,
ad_{X_{2}} \in End({\mathcal{G}}^{1}) \equiv Mat_{3}(\bf R)$.

It is obvious that $ad_{X_{1}}$ and $ad_{X_{2}}$ cannot be
concurrently trivial because ${\mathcal{G}}^{1} \cong {\bf R}^{3}$.
There is no loss of generality in assuming $ad_{X_{2}} \neq 0$. By
changing basis, if necessary, we get the similar classification of
$ad_{{X}_{2}}$ as follows:

$\begin{pmatrix} {{\lambda}_1}&0&0\\
0&{{\lambda}_2}&0\\0&0&1 \end{pmatrix},\, ({\lambda}_1, {\lambda}_2
\in {\bf R}\setminus \lbrace 1\rbrace, \, {\lambda}_1 \neq
{\lambda}_2 \neq 0$);\qquad $\begin{pmatrix}
1&0&0\\0&1&0\\0&0&{\lambda}
\end{pmatrix},\linebreak ({\lambda} \in {\bf R}\setminus \lbrace 0, 1
\rbrace$);\quad $\begin{pmatrix} {\lambda}&0&0\\0&1&0\\0&0&1
\end{pmatrix},\, ({\lambda} \in {\bf R}\setminus \lbrace 1
\rbrace$);\quad $\begin{pmatrix} 1&0&0\\0&1&0\\0&0&1
\end{pmatrix}$;\quad $\begin{pmatrix} {\lambda}&0&0\\0&1&1\\0&0&1
\end{pmatrix},\linebreak ({\lambda} \in {\bf R}\setminus \lbrace 1
\rbrace$);\qquad $\begin{pmatrix} 1&1&0\\0&1&0\\0&0&{\lambda}
\end{pmatrix},\quad ({\lambda} \in {\bf R}\setminus \lbrace 0, 1
\rbrace$);\qquad $\begin{pmatrix} 1&1&0\\0&1&1\\0&0&1
\end{pmatrix}$; \linebreak$\begin{pmatrix} cos{\varphi}&-sin{\varphi}&0\\
sin{\varphi}&cos{\varphi}&0\\0&0&\lambda \end{pmatrix},\quad
(\lambda \in {\bf R}\setminus \lbrace 0\rbrace, \, \varphi \in (0,
\pi))$. \vskip5mm

Assume that $[X_{1}, X_{2}] = mX_{3} + nX_{4} + pX_{5}; m, n, p \in
{\bf R}$. We can always change basis in order to have $[X_{1},
X_{2}] = mX_{3}$. Indeed, if
 $$ad_{X_{2}} = \begin{pmatrix} {{\lambda}_1}&0&0\\
0&{{\lambda}_2}&0\\0&0&1 \end{pmatrix},\, ({\lambda}_1, {\lambda}_2
\in {\bf R}\setminus \lbrace 1\rbrace, \, {\lambda}_1 \neq
{\lambda}_2 \neq 0),$$ then by changing $X_{1}$ for ${X_{1}}^{'} =
X_{1} + \frac{n}{\lambda_{2}}X_{4} + pX_{5}$ we get $[{X_{1}}^{'},
X_{2}] = mX_{3}$, $m \in {\bf R}$. For the other values of
$ad_{X_{2}}$, we also change basis in the same way. Hence, without
restriction of generality, we can assume right from the start that
$[X_{1}, X_{2}] = mX_{3}$, $m \in {\bf R}$.

There are three cases which contradict each other as follows.
\begin{description}
    \item[(1)] $[X_{1}, X_{2}] = 0$ ( i.e. $m = 0$ ) and
    $ad_{X_{1}} = 0$. Then $\mathcal{G} =
\mathcal{H} \oplus {\bf R}.X_{1}$, where $\mathcal{H}$ is the
subalgebra of $\mathcal{G}$ generated by $\{X_{2}, X_{3}, X_{4},
X_{5}\}$, i.e. $\mathcal{G}$ is decomposable.
    \item[(2)] $[X_{1}, X_{2}] = 0$ and $ad_{X_{1}} \neq 0$.
        \begin{description}
            \item[(2a)] Assume $ad_{X_{2}} = \begin{pmatrix} {{\lambda}_1}&0&0\\
0&{{\lambda}_2}&0\\0&0&1 \end{pmatrix},\, ({\lambda}_1, {\lambda}_2
\in {\bf R}\setminus \lbrace 1\rbrace, \, {\lambda}_1 \neq
{\lambda}_2 \neq 0$).

In view of Lemma 2.2, it follows by a direct computation that
$$ad_{X_{1}} = \begin{pmatrix} \mu&0&0\\0&\nu&0 \\ 0&0&\xi
\end{pmatrix};\, \mu, \nu, \xi \in {\bf R};\, {\mu}^{2} + {\nu}^{2} +{\xi}^{2} \ne
0.$$

If $\xi \neq 0$, by changing ${X_{1}}^{'} = X_{1} - \xi X_{2}$, we
get $$ad_{{X_{1}}^{'}} =  \begin{pmatrix} {\mu}^{'}&0&0\\0&{\nu}^{'}&0 \\
0&0&0 \end{pmatrix};$$ where \, ${\mu}^{'} = \mu - \xi
{\lambda}_{1}, {\nu}^{'} = \nu - \xi {\lambda}_{2}.$ Thus, we can
assume from the outset that $$ad_{X_{1}} = \begin{pmatrix}
\mu&0&0\\0&\nu&0
\\ 0&0&0
\end{pmatrix};\, \mu, \nu  \in {\bf R};\, {\mu}^{2} + {\nu}^{2} \ne
0.$$

Let $F = \alpha{X_{1}}^{*} + \beta{X_{2}}^{*} + \gamma{X_{3}}^{*} +
\delta{X_{4}}^{*} + \sigma{X_{5}}^{*} \in {\mathcal{G}}^{*}$ and $U
= aX_{1} + bX_{2} + cX_{3} + dX_{4} + fX_{5} \in \mathcal{G}$;\,
$\alpha, \beta, \gamma, \delta, \sigma, a, b, c, d, f \in {\bf R}$.
So we have
\begin{align} \notag {\mathcal{G}}_{F}& = KerB_{F}\\
\notag & = \{U\in \mathcal{G}/ \langle F, [U, X_{i}] \rangle = 0;\,
i = 1, 2, ,3, 4, 5\}.\end{align} Upon simple computation, we get
$$U\in {\mathcal{G}}_{F}\, \Leftrightarrow \, M
\begin{pmatrix} a\\b\\c\\d\\f\end{pmatrix} =
\begin{pmatrix}0\\0\\0\\0\\0\end{pmatrix},$$
\vskip5mm where
$$M := \begin{pmatrix}
0&0&\mu\gamma&\nu\delta&0\\
0&0&-{\lambda}_{1}\gamma&-{\lambda}_{2}\delta&-\sigma\\
\mu\gamma&{\lambda}_{1}\gamma&0&0&0\\
\nu\delta&{\lambda}_{2}\delta&0&0&0\\0&\sigma&0&0&0\end{pmatrix}.$$
\vskip5mm Hence, $dim{\Omega}_{F} = dim{\mathcal{G}} -
dim{\mathcal{G}}_{F} = rank(M)$. According to Lemma 2.3,
${\Omega}_{F}$ is a K-orbit of maximal dimension if
${F|}_{{\mathcal{G}}^{1}} \neq 0 $, i.e. if ${\gamma}^{2} +
{\delta}^{2} + {\sigma}^{2} \neq 0.$ In particular, $rank(M)$ is a
constant if $\gamma, \delta, \sigma$ are not concurrently zeros.
However, it is easily seen that $rank(M) = 2$ when $\gamma = \delta
= 0 \neq \sigma$, but $rank(M) = 4$ when all of $\gamma, \delta,
\sigma$ are different zero. This contradiction show that Case (2a)
cannot happen.
            \item[(2b)] In exactly the same way, but replacing the
            considered value of $ad_{X_{2}}$ with the others, we can
            be seen that Case (2) cannot happen anyway.
        \end{description}
    \item[(3)] $[X_{1}, X_{2}] \neq 0$ ( i.e. $m \neq 0$ ).
    By changing $X_{1}$ for ${X_{1}}^{'} = {\frac{1}{m}}X_{1}$ one has
$[{X_{1}}^{'}, X_{2}] = X_{3}$. Hence, there is no loss of
generality in assuming from the outset that $[X_{1}, X_{2}] =
X_{3}$.

By an argument similar to the one in Case (2a), we get a
contradiction again if $ad_{X_{1}} \neq 0$. In other words,
$ad_{X_{1}} = 0$. Therefore, in the dependence on the value of
$ad_{X_{2}}$, $\mathcal{G}$ will be isomorphic to one of algebras \,
${\mathcal{G}}_{5,3,1({\lambda}_{1}, {\lambda}_{2})},
({\lambda}_{1}, {\lambda}_{2} \in {\bf R}\setminus \{ 0, 1\},
{\lambda}_{1} \neq {\lambda}_{2} \neq 0 );\,
{\mathcal{G}}_{5,3,2(\lambda)}, (\lambda \in {\bf R}\setminus \{0,
1\});\, {\mathcal{G}}_{5,3,3(\lambda)}, (\lambda \in {\bf
R}\setminus \{1\});\, {\mathcal{G}}_{5,3,4};\,
{\mathcal{G}}_{5,3,5(\lambda)},\, (\lambda \in {\bf R}\setminus
\{1\});\, {\mathcal{G}}_{5,3,6(\lambda)},$\linebreak $(\lambda \in
{\bf R}\setminus \{0, 1\});\, {\mathcal{G}}_{5,3,7};\,
{\mathcal{G}}_{5,3,8(\lambda, \varphi)}, (\lambda \in {\bf
R}\setminus \{0\}, \varphi \in (0, \pi)))$. Obviously, these
algebras are not isomorphic to each  other.

To complete the proof, it remains to show that all of these algebras
are MD5-algebras. At first, we shall verify this assertion for
$\mathcal{G} = {\mathcal{G}}_{5,3,1({\lambda}_{1}, {\lambda}_{2})},
({\lambda}_{1}, {\lambda}_{2} \in {\bf R}\setminus \{ 0, 1\},
{\lambda}_{1} \neq {\lambda}_{2} \neq 0 )$. Consider an arbitrary
linear form $F = \alpha{X_{1}}^{*} + \beta{X_{2}}^{*} +
\gamma{X_{3}}^{*} + \delta{X_{4}}^{*} + \sigma{X_{5}}^{*} \in
{\mathcal{G}}^{*}$; ( $\alpha, \beta, \gamma, \delta, \sigma  \in
\linebreak {\bf R}$ ). We need prove that $dim{\Omega}_{F} =
dim{\mathcal{G}} - dim{\mathcal{G}}_{F}$ is zero or maximal.

Let $U = aX_{1} + bX_{2} + cX_{3} + dX_{4} + fX_{5} \in
\mathcal{G}$; ( $a, b, c, d, f \in {\bf R}$ ). Upon simple
computation which is similar to one in Case (2a), we get

$$U\in {\mathcal{G}}_{F}\, \Leftrightarrow \, N
\begin{pmatrix} a\\b\\c\\d\\f\end{pmatrix} =
\begin{pmatrix}0\\0\\0\\0\\0\end{pmatrix},$$
\vskip5mm where
$$N := \begin{pmatrix}
0&-\gamma&0&0&0\\
\gamma&0&-{\lambda}_{1}\gamma&-{\lambda}_{2}\delta&-\sigma\\
0&{\lambda}_{1}\gamma&0&0&0\\
0&{\lambda}_{2}\delta&0&0&0\\0&\sigma&0&0&0\end{pmatrix}.$$
\vskip5mm Hence, $dim{\Omega}_{F} = dim{\mathcal{G}} -
dim{\mathcal{G}}_{F} = rank(N)$. It is plain that

\begin{equation} rank(N) = \begin{cases} 0 & if \quad \gamma = \delta =
\sigma = 0;\\ 2 & if \quad {\gamma}^{2} + {\delta}^{2} +
{\sigma}^{2} \neq 0.\end{cases} \notag \end{equation} \vskip4mm

Therefor, ${\Omega}_{F}$ is the orbit of dimension zero or two
(maximal dimension) for any $F \in {\mathcal{G}}^{*}$, i.e.
$\mathcal{G} = {\mathcal{G}}_{5,3,1({\lambda}_{1}, {\lambda}_{2})}$
is an MD5-algebra, $({\lambda}_{1}, {\lambda}_{2} \in {\bf
R}\setminus \{ 0, 1\}, {\lambda}_{1} \neq {\lambda}_{2} \neq 0 )$.
By the same way, we can be also seen that the remaining algebras are
MD5-algebras. The proof is complete. \hfill $\square$
\end{description}

\subsection*{Concluding Remark} Let us recall that each real Lie
algebra $\mathcal{G}$ define only one connected and simply connected
Lie group G such Lie(G) = $\mathcal{G}$. Therefore we obtain a
collection of eight families of connected and simply connected
MD5-groups corresponding to given MD5-algebras in Theorem 2.1. For
convenience, each MD5-group from this collection is also denoted by
the same indices as corresponding MD5-algebra. For example,
$G_{5,3,1({\lambda}_{1}, {\lambda}_{2})}$ is the connected and
simply connected MD5-group corresponding to
${\mathcal{G}}_{5,3,1({\lambda}_{1}, {\lambda}_{2})}$. Specifically,
we have eight families of MD5-groups as follows:
$G_{5,3,1({\lambda}_{1}, {\lambda}_{2})}, ({\lambda}_{1},
{\lambda}_{2} \in {\bf R}\setminus \{ 0, 1\}, {\lambda}_{1} \neq
{\lambda}_{2} \neq 0 );\, G_{5,3,2(\lambda)}, (\lambda \in {\bf
R}\setminus \{0, 1\});\, G_{5,3,3(\lambda)}, (\lambda \in {\bf
R}\setminus \{1\});\, G_{5,3,4};\, G_{5,3,5(\lambda)}, (\lambda \in
{\bf R}\setminus \{1\}); G_{5,3,6(\lambda)}, (\lambda \in {\bf
R}\setminus \{0, 1\});\, G_{5,3,7};\, G_{5,3,8(\lambda, \varphi)},
(\lambda \in {\bf R}\setminus \{0\}, \varphi \in (0, \pi)))$. All of
them are indecomposable MD5-groups. In the next paper, we shall
describe the geometry of K-orbits of each considered MD5-group,
topologically classify MD5-foliations associated to these MD5-groups
and give a characterization of the Connes' $C^{*}$-algebras (see
[1]) corresponding to these MD5-foliations.
\subsection*{Acknowledgement}
    The author would like take this opportunity to thank their teacher,
Prof. DSc. Do Ngoc Diep for his excellent advice and support. Thanks
are due also to the first author's colleague, Prof. Nguyen Van Sanh
for his encouragement.

\end{document}